\documentclass[12pt]{amsart}

\usepackage{euscript,amssymb,xspace}
\begin{document}
\theoremstyle{plain}{
\newtheorem{theorem}{Theorem}[section]}
\theoremstyle{plain}{
\newtheorem{proposition}[theorem]{Proposition}}
\theoremstyle{plain}{
   \newtheorem{lemma}[theorem]{Lemma}}
\theoremstyle{plain}{
   \newtheorem{corollary}[theorem]{Corollary}}
\theoremstyle{definition}{
   \newtheorem{definition}[theorem]{Definition}}
\theoremstyle{definition}{
   \newtheorem{example}[theorem]{Example}}
\theoremstyle{remark}{
   \newtheorem{remark}[theorem]{Remark}}

\numberwithin{equation}{section}

\newcommand{\sst}{\scriptstyle}
\newcommand{\n}{\noindent}
\newcommand{\ds}{\displaystyle}
\newcommand{\intl}{\int\limits}
\newcommand{\vp}{\varepsilon}
\newcommand{\bb}[1]{\mathbb{#1}}
\newcommand{\cl}[1]{\mathcal{#1}}
\newcommand{\brel}[1]{\overset{\rightharpoonup}{#1}}
\newcommand{\ovl}{\overline}

\newcommand{\wdist}{\tilde{\rho}_{D}}
\newcommand{\bR}{\bb{R}}
\newcommand{\bS}{\bb{S}}

\title[$\alpha$-Continuity Properties ]{$\alpha$-Continuity Properties \\of the Symmetric $\alpha$-Stable Process}
\author{R. Dante DeBlassie} \author{Pedro J. M\'endez-Hern\'andez}

\address{Department of Mathematics, Texas A\& M University, College Station, TX 77843}
\email[Dante DeBlassie]{deblass@math.tamu.edu}
\address{Mathematics Department, University of Utah, Salt Lake City, UT 84112}
\email[Pedro M\'endez]{mendez@math.utah.edu}

\begin{abstract}
Let $D$ be a domain of finite Lebesgue measure in $\bR^d$ and let $X^D_t$ be the
symmetric $\alpha$-stable process killed upon exiting $D$. Each
element of the set $\{ \lambda_i^\alpha\}_{i=1}^\infty$ of
eigenvalues associated to $X^D_t$, regarded as a function of
$\alpha\in(0,2)$, is right continuous. In addition, if $D$ is
Lipschitz and bounded, then each $ \lambda_i^\alpha$ is continuous in $\alpha$
and the set of associated eigenfunctions is precompact. We also
prove that if $D$ is a domain of finite Lebesgue measure,  then
for all $0<\alpha<\beta\leq 2$ and $i\geq 1$,
\[\lambda_i^\alpha \leq \left[\, \lambda^\beta_i\right]^{\alpha/\beta}.\]
Previously, this bound had been known only for $\beta=2$ and $\alpha$ rational.

\end{abstract}
\maketitle

\section{Introduction}

Let $X_t$ be a $d$-dimensional symmetric $\alpha$-stable process
of order $\alpha \in (0,2]$. The process $X_t$ has stationary
independent increments and its transition density $p^{\alpha}(t,
z,w)=f_t^{\alpha}(z-w)$ is determined by its Fourier transform
\[ \exp(-t|z|^\alpha) = \int_{\bR^{d}} e^{i z\cdot
w}f_t^{\alpha}(w)dw.\]  These processes have right continuous
sample paths and their transition densities satisfy the  scaling
property
\[ p^{\alpha}(t,x,y)= t^{-d/\alpha}\,
p^{\alpha}(1,t^{-1/\alpha}x,t^{-1/\alpha} y)\;.\] When $\alpha=2$
the process $X_t$ is a $d$-dimensional Brownian motion running at
twice the usual speed. The non-local  operator associated to $X_t$
is $(-\Delta)^{\alpha/2}$ where $\Delta$ is the Laplace operator
in $\bR^d$.

Let $D$ be a domain in $\bR^d$ and let $X_t^D$ be  the symmetric
$\alpha$-stable process killed upon leaving $D$. We write
$p^{\alpha}_D(t,x,y)$ for the transition density of $X^D_t$ and
$H_\alpha$ for  its  associated  non-local self-adjoint positive
operator. It is well known that if $D$ has finite Lebesgue measure
then the spectrum of $H_\alpha$ is discrete. Let \[0<
\lambda^\alpha_1(D)< \lambda^\alpha_2(D) \leq \lambda^\alpha_3(D)
\leq \cdots\] be the eigenvalues  of $H_\alpha$, and let
\[\varphi_1^\alpha,\enskip\varphi_2^\alpha,\enskip\varphi_3^\alpha,\enskip\ldots\] be the
corresponding sequence of eigenfunctions.

Several authors have studied properties of the eigenvalues and
eigenfunctions of $H_{\alpha}$. One common theme has been to
extend results on Brownian motion ($\alpha=2$) to analogous
results for  symmetric $\alpha$-stable processes. For example, R.
M. Blumenthal and R. K. Getoor \cite{BG} have shown  Weyl's
asymptotic law holds: if $D$ is a bounded open set and
$N(\lambda)$ is the number of eigenvalues less than or equal to
$\lambda$, then there exists a constant $C_d$, depending only on
$d$, such that
\[N(\lambda)\approx C_d \frac{m(D)}{\Gamma(d+1)}\lambda^d\quad\text{ as }\lambda\to\infty,\]
provided $m(\partial D)=0$, where $m$ is Lebesgue measure.

If $D\subseteq\mathbb{R}^d$ is a domain, define the \emph{inner
radius} $R_D$ to be the supremum of the radii of all balls
contained in $D$. R. Ba\~nuelos et al \cite{BLM} and P.
M\'endez-Hern\'andez \cite{M} have shown  if $D$ is a convex
domain with finite inner radius $R_D$ and $I_D$ is the interval
$(-R_D,R_D)$, then
\[\lambda_1^\alpha(I_D)\leq\lambda_1^\alpha(D).\]

Moreover, if $D\subseteq\mathbb{R}^d$ has finite volume and $D^*$
is a ball in $\mathbb{R}^d$ with the same volume as $D$, then it
was proved in \cite{BLM} that the Faber-Krahn inequality holds:
\[\lambda_1^\alpha(D^*)\leq\lambda_1^\alpha(D).\]

Another line of inquiry taken by those authors was to consider the
eigenvalues as a function of the index $\alpha$. For instance, if
$D$ is a convex domain with finite inner radius $R_D$, then
\[\frac{2^\alpha\Gamma(1+\frac\alpha 2)\Gamma(\frac{1+\alpha}2)}{\Gamma (\frac 12)R^\alpha_D}
\leq\lambda^\alpha_1(D)\leq\lambda^\alpha_1(B_{R_D}),\] where
$B_{R_D}$ is a ball in $\bR^d$ of radius $R_D$. They also proved
that if $D\subseteq\mathbb{R}^d$ has finite volume, then
\begin{equation}\label{bnd1}\lambda_1^\alpha(D)\,\leq\,\left[ \,\mu_1(D)\,\right]^{\alpha/2},\end{equation}
where $\mu_1(D)$ is the first Dirichlet eigenvalue of $-\Delta$ on
$D$.

For the Cauchy process, i.e. $\alpha=1$,  and bounded Lipshitz
domains, R. Ba\~nuelos and T. Kulczycki \cite{Banuelos} extended
(\ref{bnd1})  to
\begin{equation}\label{bnd2}\lambda_i^1(D)\,\leq \, \left[\,\mu_i(D)\,\right]^{1/2}, \quad i=1,2,\ldots,\end{equation}
where $$0< \mu_1(D)<\mu_2(D)\leq\cdots$$ are all the Dirichlet
eigenvalues of $-\Delta$ on $D$. Their  proof of (\ref{bnd2}) is
based on a variational formula for $\lambda_i^1(D)$ they developed
from a connection with the Steklov problem for the Laplacian. They
also obtained  many detailed properties of the eigenfunctions
$\varphi_i^1$ for the Cauchy process.

By finding a connection with the symmetric stable process with
rational index $\alpha$ and  PDEs of order higher than 2, R. D.
DeBlassie \cite{DeBlassie} derived a variational formula for the
eigenvalues which led to the following extension of (\ref{bnd1})
and (\ref{bnd2}):
\begin{equation}\label{bnd3}\lambda_i^\alpha(D)\leq\left[\,\mu_i(D)\,\right]^{\alpha/2},\quad i=1,2,\ldots,\end{equation}
for all rational $\alpha\in(0,2)$ and certain bounded domains
$D\subseteq\mathbb{R}^d$. The class of admissible domains includes
convex polyhedra, Lipschitz domains with sufficiently small
Lipschitz constant and $C^1$ domains.

In this article, we study the eigenvalues and eigenfunctions
regarded as functions of the index $\alpha$. Our first result
concerns continuity of the eigenvalues.

\begin{theorem}\label{thm1.1}
Let $D$ be a domain of finite Lebesgue measure. Then, as a function of $\alpha\in
(0,2)$, $\lambda^{\alpha}_i$ is right continuous for each positive
integer $i$.
\end{theorem}

In order to prove Theorem \ref{thm1.1} we need the following interesting monotonicity
property extending (\ref{bnd3}) above.

\begin{theorem}\label{thm1.2}
Let $D$ be a domain of finite Lebesgue measure in $\mathbb{R}^n$. If $0<\alpha<\beta\leq 2$, then for all
positive integers $i$,
\[\big[\lambda^\alpha_i(D)\big]^{1/\alpha}\, \leq \,
\left[\lambda^\beta_i(D)\right]^{1/\beta}.\]
\end{theorem}

By requiring more regularity of $\partial D$, we can prove the following
extension of Theorem \ref{thm1.1}.

\begin{theorem}\label{thma}
Let $D$ be a bounded Lipschitz domain. Then, as a function of
$\alpha\in (0,2)$, $\lambda^{\alpha}_i$ is continuous for each
positive integer $i$.
\end{theorem}

We will obtain Theorem \ref{thma} from the following result that we
believe is of independent interest.

\begin{theorem}\label{thmb}
Let $D$ be a bounded Lipschitz domain. If $\alpha_m$ converges to
$\alpha\in (0,2)$ then for each positive integer  $i$,
$\{\varphi^{\alpha_m}_i\colon
\
m\ge 1\}$ is precompact in $C(\ovl D)$ equipped with the sup norm.
Moreover, if $\lambda^{\alpha_m}_i$ converges to $\lambda$, then
any limit point of $\{\varphi^{\alpha_m}_i\colon \ m\ge 1\}$ is an
eigenfunction of $H_\alpha$ and $\lambda$ is the corresponding
eigenvalue.
\end{theorem}

As a corollary of the proof of the last theorem, we obtain
continuity of the first eigenfunction as a function of $\alpha$.

\begin{theorem}\label{thm1.5}
If $D$ is a bounded Lipshitz domain and $\alpha_m$ converges to
$\alpha\in (0,2)$, then $\varphi^{\alpha_m}_1$ converges uniformly
to $\varphi^{\alpha}_1$ on $D$.
\end{theorem}

The article is organized as follows. In section 2 we present some
results needed in the proof of Theorem 1.1. We establish Theorem
1.1 in section 3 by proving upper semicontinuity and right lower
semicontinuity of the eigenvalues via Dirichlet forms. In section
4 we prove Theorem 1.2 using an extension of an operator inequality
from \cite{DeBlassie}. Lower semicontinuity of the
eigenvalues, for Lipschitz domains,  is proved in section 5 using
Theorem 1.4. This will yield Theorems 1.3 and 1.5. Section 6 deals
with certain weak convergence results needed to prove Theorem 1.4.
Finally, in section 7 we prove Theorem 1.4.

\section{Preliminary results}

Throughout this section we will assume the domain $D$ has finite
Lebesgue measure. We denote by  $C_c^{\infty}(D)$ the set of
$C^{\infty}$ functions with compact support in $D$. The inner
product and the norm in $L^2(D)$ will be denoted by
$\langle \cdot,\cdot\rangle$ and $\|\cdot\|_2$, respectively.

For any domain  $D\subseteq \bR^{d}$, we define $\tau_{D}$ to be
the first exit time of $X_t$ from $D$, i.e.,

\[ \tau_{D}=\inf \{ t>0: X_t \notin D\}.\]

Let

\[\mathcal{F}_\alpha=\left\{ \varphi \in L^2(\bR^d):\;\int\int\frac{[\,\varphi(y)-\varphi(x)\,]^2}{|y-x|^{d+\alpha}}\;dydx\,<\,\infty\,\right\}.\]

\noindent The Dirichlet form
$(\mathcal{E}_\alpha,\mathcal{F}_\alpha)$ associated to $X_t$ is
given by

\[ \mathcal{E}_{\alpha}(\psi,\varphi)=A(d,\alpha)\, \int\int\,
\frac{\left[\,\psi(y)-\psi(x)\,\right]\;\left[\,\varphi(y)-\varphi(x)\,\right]}{|y-x|^{d+\alpha}}\,dydx,\]

\noindent for all $\psi,\varphi \in \mathcal{F}_\alpha$, where

\[ A(d,\alpha)=\frac{
\Gamma(\frac{d-\alpha}{2})}{2^\alpha \,\pi^{d/2}\,\Gamma(
\frac{\alpha}{2})}.\]

\noindent It is well known that  the   Dirichlet form
corresponding to $X^D_t$ is given by $(\mathcal{E}_{\alpha},
\mathcal{F}_{\alpha,D})$, where
\[\mathcal{F}_{\alpha,D}=\left\{ u \in \mathcal{F}_\alpha:
\text{a quasi continuous version of } u \text{ is $0$ quasi
everywhere in } D^c \right\}. \]
Recall that for all $ \psi, \varphi$ in the domain of $H_\alpha$
we have

\[\mathcal{E}_{\alpha}(\psi,\varphi)=\langle \psi,H_\alpha \varphi\rangle.\]
As seen in Theorem 4.4.3 of \cite{fukushima},
$\mathcal{F}_{\alpha,D}$ is the closure of $C^\infty_c(D)$ in
$\mathcal{F}_\alpha$  with respect to the norm
\[\|\varphi\|_{\alpha}=\sqrt{\mathcal{E}_{\alpha}(\varphi,\varphi)+\|\varphi\|_2\,}.\]

\begin{lemma}
Let  $\varphi, \psi \in C^\infty_c(D)$. Then the function
\[\mathcal{E}_\alpha( \varphi,  \psi): (0,2) \to \bR,\]   is continuous on $(0,2)$.
\end{lemma}
\begin{proof}

Let $\varphi, \psi \in C^\infty_c(D)$, and   let $\beta \in
(\alpha-\delta, \alpha+\delta)$,  where $\delta= \frac{1}{2} \min
\left\{ 2-\alpha, \alpha \right\}$. Then there exists a constant
$C>0$, depending only on $\varphi$ and $\psi$, such that
\begin{eqnarray*}
\frac{\left[\,\psi(y)-\psi(x)\,\right]\;\left[\,\varphi(y)-\varphi(x)\,\right]}{|y-x|^{d+\beta}}&\leq
& \frac{C}{|y-x|^{d+\beta-2}}\\ & \leq&\,C\,\max\left\{
\frac{1}{|y-x|^{d+\alpha-\delta-2}},
\frac{1}{|y-x|^{d+\alpha+\delta-2}}\right\}.
\end{eqnarray*}
Since $D$ has finite measure a  simple computation using polar
coordinates shows
\[\max\left\{
\frac{1}{|y-x|^{d+\alpha-\delta-2}},
\frac{1}{|y-x|^{d+\alpha+\delta-2}}\right\},\] is integrable in
$D\times D$. The result immediately follows from the  dominated
convergence theorem.
\end{proof}

We end this section with some basic estimates on $L^2$ norms to be
used in the next section. Suppose $k$ is a positive integer,  $0<
\epsilon<1$, and $\varphi_1,\ldots,\varphi_k \in L^2(D)$ satisfy
\[\left|\langle \varphi_i, \varphi_j\rangle\right|\, < \frac{\epsilon}{4k^2},\qquad i\not=j, \]
\[
\left(\, 1-\frac{\epsilon}{4k^2}\,\right)\,<\,\|\varphi_i\|_{2}^2< \left(\,
1+\frac{\epsilon}{4k^2}\,\right),\]

\bigskip
\noindent for all $ 1\leq i,j\leq k$. If $ \psi= \sum_{i=1}^k a_i
\varphi_i,$ with $ \|\psi\|_2=1,$ then we now show
\begin{equation}\label{norm1}\frac{1}{1+\epsilon/2}\leq\sum_{i=1}^k
a_i^2\leq \frac{1}{1-\epsilon/2},\end{equation} and
$\varphi_1,\dots\varphi_k$ are linearly independent.

For the proof, note we have
\begin{eqnarray*}
1 =\; \langle \psi,\psi\rangle &=& \sum_{i=1}^{k}
a_i^2\|\varphi_i\|^2_2+ 2 \sum_{i=1}^k\sum_{j>i}^k a_i a_j \langle
\varphi_i, \varphi_j\rangle\\ &\geq&\sum_{i=1}^{k}  a_i^2
\,\left(\, 1 -\frac{\epsilon}{4k^2}\,\right)\,-\,2
\sum_{i=1}^k\sum_{j>i}^k |a_i|\,|a_j|\,\frac{\epsilon}{4k^2}\\
&\geq&\sum_{i=1}^{k} a_i^2 \,\left( \,1
-\frac{\epsilon}{4k^2}\,\right)- \,(k^2-k)\, \sum_{i=1}^k
a_i^2\frac{\epsilon}{4k^2}\\ &\geq&
\left(1-\epsilon/2\right)\,\sum_{i=1}^k a_i^2,
\end{eqnarray*}
and we conclude
\begin{equation*}\sum_{i=1}^k
a_i^2\leq \frac{1}{1-\epsilon/2}.\end{equation*}

\noindent Similar computations give the remaining assertions.

\section{Proof of Theorem \ref{thm1.1}}

We will use the following well known result, see \cite{davies2}.

\begin{theorem}\label{thm3.1} Let $H$ be a non-negative self-adjoint unbounded operator
with discrete spectrum $\{\lambda_i\}_{i=1}^\infty$, and domain
$Dom(H)$. Then for $i \geq 1$
\begin{equation}\lambda_i = \inf \left\{ \lambda(L):
L \subseteq Dom(H), dim(L)=i\right\},\end{equation} where
\begin{equation} \lambda(L)= \sup \left\{ \langle  Hf ,f\rangle \,:\, f \in L,
\|f\|_2=1\right\},\label{L}\end{equation} and $L$ is a vector
subspace of $Dom(H)$ of dimension $i$.
\end{theorem}

We will prove the right continuity of the $k$-th eigenvalue in several
steps.

\begin{proposition}\label{prop3.2}
Let $D$ be a domain of finite Lebesgue measure. Then for all $k
\geq 1$
\[\limsup_{\beta \to \alpha} \lambda^\beta_k(D) \leq
\lambda^\alpha_k(D).\]
\end{proposition}
\begin{proof}
Let $0< \epsilon <1$ and $k\ge 1$. Recall  $C^\infty_c(D)$ is
dense in $Dom(H_\alpha)$ under the norm $\|\cdot\|_\alpha$. Then
for all $\alpha \in (0,2)$, there exist $\varphi_1, \ldots,
\varphi_k \in C^\infty_c(D)$ such that
\begin{equation}
|\langle \varphi_i^\alpha,\varphi_j^\alpha\rangle-\langle\varphi_i,\varphi_j\rangle| \,<\,
\frac{\epsilon}{8k^2},\label{equa1}\end{equation} and
\begin{equation}|\mathcal{E}_\alpha(\varphi_i^\alpha,\varphi_j^\alpha )-\mathcal{E}_\alpha
(\varphi_i,\varphi_j)|\,<\,
\frac{\epsilon}{8k^2},\label{equa2}\end{equation} for all $1\leq
i,j\leq k$.

Thanks to Lemma 2.1 there exists  $\eta_0$  such that for all
$\beta \in (\alpha-\eta_0, \alpha+\eta_0)$
\begin{equation}|\mathcal{E}_{\alpha}(\varphi_i,\varphi_j) \,-\,
\mathcal{E}_{\beta}(\varphi_i,\varphi_j)|\,<\,
\frac{\epsilon}{8k^2}.\label{equa3}\end{equation}

Notice  (\ref{equa1}) implies
\[|\langle\varphi_i,\varphi_j\rangle| \,<\,
\frac{\epsilon}{8k^2},\quad i\not= j,\] and
\[1-\frac{\epsilon}{8k^2}\,<\,\|\varphi_i\|_2^2 \,<\,
1+\frac{\epsilon}{8k^2},\] for all $1\leq i,j \leq k$. Then by the
comments at the end of section 2, we know
$\varphi_1,\dots,\varphi_k$ are linearly independent.

Theorem \ref{thm3.1}  implies
\[ \lambda^\beta_k(D)\leq \lambda_\beta(L_k),\] where $L_k=\text{span}\{\varphi_1,\ldots,\varphi_k\}$ and
\[\lambda_\beta(L_k)=\sup\left\{\,\langle H_\beta f,f\rangle :f\in L_k,\|f\|_2=1\,\right\}.\]

Take  $\psi= \sum_{i=1}^k a_i \varphi_i \in L_k$  such that

\begin{equation}\label{equa4}\lambda_\beta(L_k) \leq \mathcal{E}_\beta(\psi,\psi)+
\epsilon/4,\end{equation} and
\[\|\psi\|_2=1.\]

Thanks to (\ref{norm1}), with $\epsilon$ there replaced by
$\epsilon/2$, we have
\[ \sum_{i=1}^k a_i^2\leq 2.\]

Then since
\[ |\,\mathcal{E}_\beta(\psi,\psi)-
\mathcal{E}_\alpha(\psi,\psi)\,| \leq
\sum_{i=1}^k\sum_{j=1}^k\,|a_i\,a_j| \;\left|\,
\mathcal{E}_\beta(\varphi_i,\varphi_j)-\mathcal{E}_\alpha(\varphi_i,\varphi_j)\,\right|,
\] (\ref{equa3}) implies
\begin{equation} \label{equa6}|\,\mathcal{E}_\beta(\psi,\psi)- \mathcal{E}_\alpha(\psi,\psi)\,| <
\frac{\epsilon}{4}.\end{equation} Thus

\begin{equation} \label{equa7}\lambda^\beta_k(D)\leq
\mathcal{E}_\alpha(\psi,\psi)+\epsilon/2.\end{equation}

Consider $\psi_0 =\sum_{i=1}^k a_i \varphi_i^\alpha$. By
(\ref{norm1}) we have
\[ \frac{1}{1+\epsilon/4}\leq  \|\psi_0\|^2_2= \sum_{i=1}^k a_i^2 \leq\frac{1}{1-\epsilon/4}
.\]

Following the argument used to obtain (\ref{equa6}), one easily
proves (\ref{equa2}) implies
\[ |\,\mathcal{E}_\alpha(\psi_0,\psi_0)-
\mathcal{E}_\alpha(\psi,\psi)\,| < \frac{\epsilon}{4}.\]
Hence
\begin{eqnarray*}
\lambda^\beta_k(D)&\leq&\mathcal{E}_\alpha(\psi,\psi)+\epsilon/2\\
&\leq&\mathcal{E}_\alpha(\psi_0,\psi_0)+3\epsilon/4\\
&=&\sum_{i=1}^ka_i^2\lambda^\alpha_i(D)+3\epsilon/4\\
&\leq&\lambda^\alpha_k(D)\sum_{i=1}^ka_i^2+3\epsilon/4\\
&\leq&\frac{1}{1-\epsilon/4}\lambda^\alpha_k(D)+3\epsilon/4,
\end{eqnarray*}
and the result immediately follows.
\end{proof}

\begin{proposition}\label{prop3.3} Let $D$ be a domain of finite Lebesgue
measure. If $D_n\subset\subset D$ is a sequence of $C^\infty$ domains
increasing to $D$, then for all $k\geq 1$,
\[\lim_{n\to\infty}\lambda_k^\alpha(D_n)=\lambda_k^\alpha(D).\]
\end{proposition}
\begin{proof}
By domain monotonicity,
\[\lambda_k^\alpha(D)\leq\lambda_k^\alpha(D_n).\]
Hence it suffices to show
\begin{equation}\label{eq4a}\limsup_{n\to\infty}\lambda_k^\alpha(D_n)\leq\lambda_k^\alpha(D).\end{equation}
To this end, let $\eta
>0$. Following the arguments presented above, we can prove
there exist $k$ linearly independent functions
$\varphi_1,\ldots,\varphi_k \in C^\infty_c(D)$ such that

\[\lambda_k^\alpha(D) \geq \lambda_\alpha(L_k)-\eta ,\]
where $L_k$ is the vector space generated by $\{\varphi_1,\ldots,\varphi_k\}$
and
\[\lambda_\alpha(L_k)=\sup\,\left\{\langle H_\alpha f,f\rangle:f\in L_k,\|f\|_2=1\,\right\}.\]

Then there exists $n_0$ such that the supports of
$\varphi_1,\ldots,\varphi_k$ are contained in $D_{n_0}$.
Consequently for large $n$, Theorem \ref{thm3.1} implies

\[\lambda_k^\alpha(D) \geq \lambda_\alpha(L_k)-\eta\geq
\lambda_k^\alpha(D_n)-\eta.\]
Hence upon letting $n\to\infty$ and $\eta\to 0$, we get (\ref{eq4a}).
\end{proof}

\begin{proposition}\label{prop3.4}
Let $D$ be a domain of finite Lebesgue measure. Then for all $k\ge 1$
\[\liminf_{\beta \to \alpha+} \lambda^\beta_k(D) \geq
\lambda^\alpha_k(D).\]
\end{proposition}
\begin{proof}
%Let us first assume the diameter of $D$ is less than 1. Then for all
%$x,y \in D$, and all $\epsilon>0$  we have
%
%\[\frac{1}{|x-y|^{d+\alpha}}\leq
%\frac{1}{|x-y|^{d+\alpha+\epsilon}}.\]
%%
%It follows that
%\[
%\int_D\int_D\frac{[\,\varphi(y)-\varphi(x)\,]^2}{|y-x|^{d+\alpha}}dxdy\leq
%\int_D\int_D\frac{[\,\varphi(y)-\varphi(x)\,]^2}{|y-x|^{d+\alpha+\epsilon}}dxdy,\]
%for all $\varphi \in \mathcal{F}_{\alpha+\epsilon,D}$.

%On the other hand Theorem 1 in \cite{adams}  states that any set
%of zero $(\alpha+\epsilon)$-Riesz capacity also has zero
%$\alpha$-Riesz capacity. Thus
%$\mathcal{F}_{\alpha+\epsilon,D}\subseteq\mathcal{F}_{\alpha,D}$.
%We conclude
%\[ \lambda_\alpha(L)\leq \frac { A(d,\alpha)}{A(d,\alpha+\epsilon)}\; \lambda_{\alpha+\epsilon}(L)\]
%for all $k$-dimensional vector subspaces $L$ of
%$\mathcal{F}_{\alpha+\epsilon,D}$. It immediately follows that

%begin{equation} \lambda^\alpha_k(D)\, \leq\, \frac {
%A(d,\alpha)}{A(d,\alpha+\epsilon)}
%\;\lambda_k^{\alpha+\epsilon}(D),\label{eq3.9}\end{equation}

%For general bounded $D$,   let $R$ be the diameter of $D$. Then
%scaling and (\ref{eq3.9}) imply

%\begin{equation} \lambda^\alpha_k(D)\,\leq \,\frac{A(d,\alpha)R^{\epsilon}}
%{ A(d,\alpha+\epsilon)}\;
%\lambda_k^{\alpha+\epsilon}(D),\label{eq3.10}\end{equation} the
%result immediately follows after letting $\epsilon \to 0$.

Let  $D$  be a domain with  finite Lebesgue measure, and let $D_n\subset\subset D$ be a
sequence of bounded $C^\infty$ domains increasing to $D$. Such a sequence can be
constructed using the regularized distance function---see page 171
in \cite{stein}.

By Theorem \ref{thm1.2},
\[\lambda_k^\alpha(D_n)\leq\left[\lambda_k^{\alpha+\epsilon}(D_n)\right]^{\alpha/(\alpha+\epsilon)}.\]
There is no danger of circular reasoning here because the proof of Theorem \ref{thm1.2} given in
the next section is independent of Theorem \ref{thm1.1}.
Now let $n\to\infty$ and appeal to Proposition \ref{prop3.3} to get
\[\lambda_k^\alpha(D)\leq\left[\lambda_k^{\alpha+\epsilon}(D)\right]^{\alpha/(\alpha+\epsilon)}.\]
Upon letting $\epsilon\to 0$, we get the desired $\liminf$ behavior.

\end{proof}

Combining Propositions \ref{prop3.2} and \ref{prop3.4}, we get Theorem \ref{thm1.1}.

\section{Proof of Theorem \ref{thm1.2}}

The first result we need is the following extension of Theorem 1.3 in
\cite{DeBlassie}. It says the operator $e^{-(-\Delta)^{\alpha/2} t}$
dominates $e^{-H_\alpha t}$ on $L^2(D)$.

\begin{theorem}\label{thm4.1}Let $D$ be a bounded smooth domain.
If $0 < \alpha \leq 2$ then for all  $\psi \in L^2(D)$,
\[\langle\psi,e^{-(-\Delta)^{\alpha/2} t}\psi\rangle\;\geq\; \langle\psi,e^{-H_\alpha t}\psi\rangle.
\]
\end{theorem}
\begin{proof}
The proof of the case given in section 6 of \cite{DeBlassie} goes through with the following changes. All
lemma and equation labels are from that article. In the proof of Lemma 6.1, it is enough to use the bound from
Lemma 3.2 in place of (5.3).

The second expression in the scaling relation (6.2) should be replaced by
\[b_\ell=a_\ell/M^\alpha\]
and any subsequent appearance of $a_\ell/M^2$ should be replaced by $a_\ell/M^\alpha$.

In the proof of Lemma 6.2 it is not necessary to appeal to the Weyl Asymptotic formula. It is
enough that $a_\ell\to\infty$ as $\ell\to\infty$ and then one needs only examine the function
$F(x)=x^ne^{-Bx}$ instead of $x^{2n/d}e^{-Bx^{2/d}}$.

The only other change is at the end of section 6 where it is shown $g(1/\delta)\to 0$ as $\delta\to 0$.
This time use the bound
\[p^\alpha(t,x,y)\,\leq\, \frac{c_{\alpha,d}}{ t^{d/\alpha}}.\]
\end{proof}

With the aid of this Theorem, we can now prove Theorem \ref{thm1.2}
for bounded smooth domains.
\begin{theorem}\label{thm4.2} Let $D$ be a bounded smooth domain in $\mathbb{R}^n$.
If $0<\alpha<\beta\leq 2$, then for all positive integers $i$,
\[\big[\lambda^\alpha_i(D)\big]^{1/\alpha}\, < \,
\left[\lambda^\beta_i(D)\right]^{1/\beta}.\]
\end{theorem}
\begin{proof}

Let $0 < \alpha< \beta\leq 2$. Denote by  $\sigma_t$  the stable subordinator of index $\frac{\alpha}{\beta}$  with
\[ E \left[ \, e^{-s \sigma_t}\,\right]= e^{-t \, s^{\alpha/\beta}},\, \text{ for all } \, s>0.\]
It is well known that for $ \psi \in L^2(D)$,
\[ E\left[\, \langle\psi,e^{-(-\Delta)^{\beta/2} \sigma_t}\psi\rangle \,\right]= \langle\psi,e^{-(-\Delta)^{\alpha/2} t}\psi\rangle,\]
and
\[ E\left[\, \langle\psi,e^{- H_\beta \sigma_t}\psi\rangle \,\right]= \langle\psi,e^{-(H_\beta)^{\alpha/\beta} t}\psi\rangle\]
(see Example 32.6  in  \cite{sato}). Then taking
\[\psi=\sum_{j=1}^{i}a_j\varphi_j^\beta \quad\text{ with }\quad \|\psi\|_2=\sum_{j=1}^ia_j^2=1,\]
Theorem \ref{thm4.1} implies
\begin{eqnarray*} \nonumber
\langle\psi,e^{-(-\Delta)^{\alpha/2} t}\psi\rangle&\geq&
\langle\psi,e^{-(H_\beta)^{\alpha/\beta} t}\psi\rangle \\
&=&\sum_{j,k=1}^ia_ja_k\,
\langle\varphi_j^\beta,e^{-(H_\beta)^{\alpha/\beta}t}\varphi_k^\beta\rangle\\
&=&\sum_{j,k=1}^ia_ja_k \, E\left[\, \langle\varphi_j^\beta,e^{-
H_\beta\sigma_t}\varphi_k^\beta\rangle \,\right]\\
&=&\sum_{j,k=1}^ia_ja_k \,E\left[\, \delta_{jk}\, e^{-\sigma_t
\lambda_j^\beta(D)}\,\right]\\ &=&  \sum_{j=1}^i a_j^2 \,e^{-t
\left[\,\lambda_j^{\beta}(D)\,\right]^{\alpha/\beta}}\\ &\geq
&\sum_{j=1}^i a_j^2\,e^{-t
\left[\,\lambda_i^{\beta}(D)\,\right]^{\alpha/\beta}}\\ &=&e^{-t
\left[\,\lambda_i^{\beta}(D)\,\right]^{\alpha/\beta}}.
\end{eqnarray*}

On the other hand, Theorem 1 in \cite{adams}  states that any set
of zero $\beta$-Riesz capacity also has zero
$\alpha$-Riesz capacity, i.e.,
$\mathcal{F}_{\beta,D}\subseteq\mathcal{F}_{\alpha,D}$.
Thus, by the formula
\[ \mathcal{E}_\alpha(\psi,\psi)= \lim_{t \to \infty}\frac{1}{t} \big\langle\psi,\left(\,1-e^{- t (-\Delta)^{\alpha/2}  }\right)\,\psi\big\rangle ,\]
we conclude
\begin{equation} \mathcal{E}_\alpha \left( \psi,\psi\right) \leq \left[\,\lambda_i^{\beta}(D)\,\right]^{\alpha/\beta},\label{monotonicity}\end{equation}
for all $ i\geq 1$.

The desired bound follows from (\ref{monotonicity}) and Theorem \ref{thm3.1}.
\end{proof}

Now we can prove Theorem \ref{thm1.2}.
Let  $D\subseteq\mathbb{R}^d$ be a domain of finite Lebesgue
measure. Suppose $D_n\subset\subset D$ is a sequence of bounded
$C^\infty$ domains increasing to $D$.

Thanks to Theorem \ref{thm4.2}, if $0<\alpha<\beta\leq 2$ then
for all positive integers $n$ and $i$,

\[ \left[\,\lambda^\alpha_i(D_n)\,\right]^{1/\alpha}\,\leq\, \left[\,
\lambda^\beta_i(D_n)\,\right]^{1/\beta}.\]
Letting $n\to\infty$, Proposition \ref{prop3.3} implies
\[ \left[\,\lambda^\alpha_i(D)\,\right]^{1/\alpha}\,\leq\, \left[\,
\lambda^\beta_i(D)\,\right]^{1/\beta},\]
as desired.\hfill$\square$

\section{Proof of Theorems 1.3 and 1.5}

\indent

We now show how Theorem \ref{thmb} implies Theorem \ref{thma}. In
order to simplify the notation, throughout this section we will
write  $\lambda^{\alpha}_k$ for  $\lambda^{\alpha}_k(D)$ and
$\mu_k$ for $\mu_k(D)$.

\begin{proof}[Proof of Theorem \ref{thma}]

We proceed by induction on $i$. For $i=1$, let
$\{\alpha_m\}_{m=1}^\infty$ be a sequence  converging to $\alpha$
in $(0,2)$. Consider any subsequence $\beta_r = \alpha_{m_r}$.
Theorem 1.2 implies the sequence
$\{\lambda^{\alpha_m}_1\}_{m=1}^\infty$ is bounded, and so there
is a subsequence $\gamma_\ell = \beta_{r_\ell}$ such that
$\lambda^{\gamma_\ell}_1$ converges as $\ell\to\infty$, say to
$\lambda$. Thanks to Theorem \ref{thmb} we can choose a
subsequence $\eta_p = \gamma_{\ell_p}$ such that
$\varphi^{\eta_p}_1$ converges uniformly to $\varphi$  an
eigenfunction of $H_\alpha$ with eigenvalue $\lambda$. Since
$\varphi^{\eta_p}_1$ is nonnegative, so is $\varphi$. But the only
nonnegative eigenfunction of $H_\alpha$ is $\varphi^{\alpha}_1$.
Thus $\lambda = \lambda^{\alpha}_1$ and $\varphi =
\varphi^{\alpha}_1$. Hence we have shown any subsequence of
$\lambda^{\alpha_m}_1$ contains a further subsequence  converging
to $\lambda^{\alpha}_1$. We conclude that $$\lim_{m \to
\infty}\lambda^{\alpha_m}_1 = \lambda^{\alpha}_1.$$ Note this also
proves Theorem 1.5.

Next assume the theorem is true for $j\le i$. We verify it is true
for $j=i+1$. We will show
\begin{equation} \label{eq5.1}\liminf\limits_{\beta\to\alpha} \lambda^{\beta}_{i+1}\ge
\lambda^{\alpha}_{i+1}.\end{equation} Combined with the lim sup
behavior from Proposition 3.2, we conclude the desired result
\[\lim_{\beta \to \alpha}\lambda^{\beta}_{i+1} = \lambda^{\alpha}_{i+1}.\]

To get the lim inf behavior, by way of contradiction, assume
$\lambda = \liminf\limits_{\beta\to \alpha} \lambda^{\beta}_{i+1}
< \lambda^{\alpha}_{i+1}$. Let $\{\alpha_m\}_{m=1}^\infty$ be a
sequence converging to $\alpha$ with
\begin{equation}
\lim_{m \to \infty}\lambda^{\alpha_m}_{i+1}=
\lambda.\label{eq5.2}\end{equation} By the induction hypothesis,
$\lambda^{\alpha_m}_j$ converges to $\lambda^{\alpha}_j$ for $j\le
i$. Then Theorem \ref{thmb} implies  we can choose a subsequence
$\beta_r = \alpha_{m_r}$ such that:

\begin{itemize}

\item For each $j$, $1\leq j\leq i$, $\lambda^{\beta_r}_j$ converges to $
\lambda^{\alpha}_j$, and $\varphi^{\beta_r}_j$ converges uniformly
to an eigenfunction $\varphi_j$ of $H_\alpha$ with corresponding
eigenvalue $\lambda^\alpha_j$.

\item The limit $\lambda$ from (\ref{eq5.2}) is an eigenvalue of $H_\alpha$, and
$\varphi^{\beta_r}_{i+1}$ converges uniformly to an eigenfunction
$\varphi_{i+1}$ of $H_\alpha$ with eigenvalue $\lambda$.

\end{itemize}

Since $\lambda$ is an eigenvalue strictly less than
$\lambda^{\alpha}_{i+1}$, we can choose positive integers $\ell$
and $m$ such that  $\ell \le m \le i$, $\lambda^{\alpha}_m =
\lambda$ and
\[
\lambda^{\alpha}_{\ell-1} < \lambda^{\alpha}_\ell =\cdots=
\lambda^{\alpha}_m <
\lambda^{\alpha}_{m+1}\leq\cdots\leq\lambda^{\alpha}_{i+1}.
\]
In particular, if $E$ is the eigenspace corresponding to $\lambda
= \lambda^{\alpha}_m$,
\[
\dim(E) = m-\ell+1.
\]

On the other hand, the uniform convergence implies for $j_1,j_2
\in \{1,\ldots, i+1\}$
\[
\delta_{j_1j_2}  = \intl_D \varphi^{\beta_r}_{j_1}
\varphi^{\beta_r}_{j_2} \ dx \,\text{ converges to }\, \intl_D
\varphi_{j_1} \varphi_{j_2}\ dx.
\]
Thus $\{\varphi_1,\ldots,\varphi_{i+1}\}$ is an orthonormal set,
and so $\{\varphi_\ell,\ldots, \varphi_m\}\cup \{\varphi_{i+1}\}$
is an orthonormal subset of $E$. This forces $\dim(E) \ge
m-\ell+2$, which contradicts (\ref{eq5.2}).  We conclude
(\ref{eq5.1}) holds.
\end{proof}

\section{Weak Convergence Results}

\indent

Let ${\bb D}[0,\infty)$ be the space of right continuous functions
$\omega\colon \ [0,\infty)\to {\bb R}^d$ with left limits. That
is, $\omega(t^+) = \lim\limits_{s\to t^+} \omega(s)  = \omega(t)$
and $\omega(t^-) = \lim\limits_{s\to t^-} \omega(s)$ exists. The
usual convention is $\omega(0^-) := \omega(0)$. Let $X_t(\omega) =
\omega(t)$ be the coordinate process and let ${\cl F}_t$ be the
$\sigma$-algebra generated by the cylindrical sets. We equip ${\bb
D}[0,\infty)$ with the Skorohod topology. Our main reference is
Chapter 3 in Ethier and Kurtz \cite{Ethier}. Let $P^{\alpha}_x$
denote the law on ${\bb D}[0,\infty)$ of the symmetric
$\alpha$-stable process started at $x$; the corresponding
expectation will be denote by $E^{\alpha}_x$.

\begin{lemma}\label{lem6.1}
If $(x_n,\alpha_n)$ converges  to $(x,\alpha)$ in ${\bb R}^n
\times (0,2)$, then $P^{\alpha_n}_{x_n}$ converges weakly to
$P^{\alpha}_x$  in ${\bb D}[0,\infty)$.
\end{lemma}

\begin{proof}
Using characteristic functions it is easy to show the
corresponding finite dimensional distributions converge. Thus by
Theorem 7.8 on page 131 in \cite{Ethier}, it suffices to show
$\{P^{\alpha_n}_{x_n}\colon \ n\ge 1\}$ is tight. For this, note
for $\beta = \alpha$ or $\alpha_n$ and $y = x$ or $x_n$,
$P^{\beta}_y$ solves the martingale problem:
\begin{itemize}
\item[a)] $P^{\beta}_y(X_0 = y) = 1$
\item[b)] for each $f\in C^2_b({\bb R}^d)$,
\end{itemize}
\[
f(X_t) - f(X_0) - \int^t_0 {\cl L}_\beta f(X_s)ds
\]
is a $P^{\beta}_y$-martingale, where $C^2_b({\bb R}^d)$ is the
space of functions with bounded continuous derivatives up to and
including order 2 and
\[
{\cl L}_\beta f(x) = A(d,\alpha) \intl_{{\bb R}^d\backslash \{x\}}
\frac{f(y) - f(x) - \nabla f(x)\cdot (y-x)
I(|y-x|<1)}{|y-x|^{d+\alpha}}dy,
\]
see section 2 of \cite{Bass2}. It is easy to show for any $f\in
C^2_b({\bb R}^d)$ there exists $C_f>0$ independent of $\alpha_n$
and $x_n$ such that $f(X_t) - f(X_0)-C_ft$ is a
$P^{\alpha_n}_{x_n}$-supermartingale. Then by Proposition 3.2 in
\cite{Bass1}, $\{P^{\alpha_n}_{x_n}\colon \ n\ge 1\}$ is tight on
${\bb D}[0,t_0]$ for all $t_0$. Even though that result is stated
in one dimension and $x_n\equiv x$, it is easy to check the proof
works in higher dimensions with $x_n$ converging to $x$.
\end{proof}

The next step is to show for each $T>0$ the distribution of
$X_{T\wedge \tau_D}$ under $P^{\alpha_n}_{x_n}$ converges to that
under $P^{\alpha}_x$ as $(x_n,\alpha_n)$ converges to
$(x,\alpha)$. To this end, define

\[{\cl A}_D = \left\{\,\omega\in {\bb D}[0,\infty)\colon \ d(X[0,\tau_D(\omega)-r], D^c)> 0
\text{ for all rational } r<\tau_D(\omega)\,\right\},
\]
and
\[
{\cl C}_D = \left\{\,\omega\in {\bb D}[0,\infty)\colon \
X(\tau_D(\omega))\in \ovl D^c\,\right\} \cap {\cl A}_D.
\]
Here $X[0,t] = \{X_s\colon \ 0\le s\le t\}$ and $d(A,B)$ is the
distance between $A$ and $B$.

\begin{lemma}\label{lem6.2}
For open $D\subseteq {\bb R}^d$, $\tau_D$ is continuous on ${\cl
C}_D$.
\end{lemma}

\begin{proof}
Let $\omega\in {\cl C}_D$ and suppose $\omega_n$ converges to $
\omega$ in ${\bb D}[0,\infty)$. We will show $\tau_D(\omega_n)$
converges to $\tau_D(\omega)$. Let
\[
\Lambda' = \left\{\,\lambda\colon \ [0,\infty)\to [0,\infty)\mid
\lambda \text{ is strictly increasing and surjective}\,\right\}
\]
Proposition 5.3 (a) and (c), on page 119, in \cite{Ethier} implies
that for each $T>0$ there exist $\lambda_n\in \Lambda'$ such that
\begin{equation}\label{eq6.1}
\lim_{n \to \infty} \sup_{0\le t\le T} |\lambda_n(t)-t|= 0,
\end{equation}
\begin{equation}\label{eq6.2}\lim_{n \to \infty} \sup_{0\le t\le T}
|\,\omega_n(t) - \omega(\,\lambda_n(t)\,)\,| = 0.
\end{equation}

First we show
\begin{equation}\label{eq6.3}
\liminf_{n\to\infty} \tau_D(\omega_n) \ge \tau_D(\omega).
\end{equation}
Let $\delta\in (0,\tau_D(\omega)/2)$ be rational and set
\[
\vp = d(\,\omega[0, \tau_D(\omega) - \delta], D^c\,).
\]
Since $\omega\in {\cl C}_D$ we have $\vp>0$. Using $T =
\tau_D(\omega)$ in (\ref{eq6.1})--(\ref{eq6.2}), there exists $N$
such that for $n\ge N$,
\begin{equation}\label{eq6.4}
\left\{\begin{array}{l} t-\delta < \lambda_n(t) < t+\delta \quad
\text{for all}\quad t\le T,\\ \noalign{\medskip} \sup\limits_{0\le
t\le T} |\,\omega_n(t) - \omega(\lambda_n(t))\,| <
\dfrac\vp2.\end{array} \right.
\end{equation}
In particular, for all $t\le \tau_D(\omega)-2\delta$ and $n\ge N$
\[
\lambda_n(t) < t+\delta \le \tau_D(\omega)-\delta<T.
\]
Thus $\omega(\lambda_n(t)) \in D$ and $d(\,\omega(\lambda_n(t))\,,
D^c\,) \ge \vp$. Therefore \[\omega_n(t) \in
B(\,\omega(\lambda_n(t))\,,\epsilon/2\,)\subseteq D,\] for all
$t\le \tau_D(\omega)-2\delta$ and $n\ge N$. This implies
$\tau_D(\omega_n) > \tau_D(\omega) -2\delta$ for $n\ge N$. Take
the lim inf as $n\to \infty$ and then let $\delta \to 0$ to get
\eqref{eq6.3}.

To finish, we show
\begin{equation}\label{eq6.5}
\limsup_{n\to\infty} \tau_D(\omega_n) \le \tau_D(\omega).
\end{equation}
Given that  $\omega\in {\cl C}_D$ and   $\omega$ is right
continuous, we can  choose $\delta>0$ such that
\[
\vp := d(\;\omega[\,\tau_D(\omega), \tau_D(\omega) + 2\delta\,],
D\;)
> 0.
\]
Using $T = \tau_D(\omega) + 2\delta$ in
(\ref{eq6.1})--(\ref{eq6.2}) we can choose $N$ such that for $n\ge
N$, \eqref{eq6.4} holds for this choice of $\delta,\vp$ and $T$.
In particular, for $n\ge N$,
\[
\tau_D(\omega)<\lambda_n(\tau_D(\omega)+\delta) <
\tau_D(\omega)+2\delta,
\]
and
\[
|\,\omega_n(\,\tau_D(\omega)+\delta\,)\, -\,
\omega(\,\lambda_n(\tau_D(\omega)+\delta)\,)\,| < \frac\vp2.
\]
Together these imply
\[
d(\omega_n(\tau_D(\omega)+\delta), D) > 0, \quad n\ge N,
\]
which in turn yields
\[
\tau_D(\omega_n) \le \tau_D(\omega)+\delta,\quad n\ge N.
\]
Taking the lim sup as $n\to\infty$ and then letting $\delta\to 0$
yields \eqref{eq6.5}.
\end{proof}

\begin{lemma}\label{lem6.3}
If $D$ is a  bounded domain that  satisfies an exterior cone
condition, or if $D$ is a cone. Then for all $x\in D$ and
$0<\alpha<2$,
\[
P^{\alpha}_x(\,{\cl C}_D\cap \{X(\tau^-_D) \in D\}\,)  = 1.
\]
\end{lemma}

\begin{proof}
If $D$ is bounded and satisfies a uniform exterior cone condition,
it is known that
\begin{equation}\label{eq6.6}
P^{\alpha}_x(\,X(\tau_D)\in \partial D\,) = 0,
\end{equation}
 see Lemma 6 in \cite{Bogdan}. If $D$ is a cone,  we can
apply   Lemma 6 in \cite{Bogdan} to $D\cap B_M(0)$ and letting
$M\to\infty$, we get \eqref{eq6.6}.

The proof of Theorem 2 in \cite{Ikeda}  implies
\begin{equation}\label{eq6.7}
P^{\alpha}_x(\,X(\tau^-_D) \in \partial D, X(\tau_D)\in E\,)  =
0,\qquad E\subseteq \ovl E\subseteq\ovl D^c,
\end{equation}
(see the lines before the footnote on page 89). Combined with
\eqref{eq6.6},
\[
P^{\alpha}_x(\,X(\tau_D)\in \ovl D^c, X(\tau^-_D)\in D\,) = 1.
\]
Thus to prove the lemma we need to show
\begin{equation}\label{eq6.8}
P^{\alpha}_x(\,d(X[0,\tau_D-r],D^c\,)>0 \text{ for all rational }
r<\tau_D\,) = 1.
\end{equation}
Let
\[
D_n = \left\{x\in D\colon \ d(x,D^c) > \frac1n\right\}
\]
and observe $\tau_{D_n}\le \tau_D$ increases to some limit $L\le
\tau_D$. By  quasi-left continuity, $X(\tau_{D_n})\to X(L)$ almost
surely. One  easily sees $X(L) \notin D$, i.e., $\tau_D\le L$.
Hence $\tau_D=L$, and the increasing limit of $\tau_{D_n}$  is
$\tau_D$.

If for some rational $r<\tau_D$ we have
\[
d(X[0,\tau_D-r], D^c) = 0
\]
then for some sequence $s_n\le \tau_D-r$,
\[
d(\,X_{s_n}, D^c\,) \to 0.
\]
It is no loss to assume $s_n$ converges, say to $s$. Choose $N$
such that for all $n\ge N$
\[
\tau_D-r < \tau_{D_n} \le \tau_D.
\]
Given $n\ge N$, choose $M_n$ such that for all $m\ge M_n$,
\[
 d(\,X_{s_m}, D^c\,) < \frac1{2n}.
\]
Then for such $m$, $X_{s_m}\in D^c_n$, which forces
\[
\tau_{D_n} \le s_m \le \tau_D-r.
\]
Let $m\to\infty$ to get $\tau_{D_n} \le s\le \tau_D-r$, then let
$n\to\infty$ to get $\tau_D = \lim\limits_{n\to\infty} \tau_{D_n}
\le \tau_D-r$; contradiction. Thus \eqref{eq6.8} holds.
\end{proof}

\begin{lemma}\label{lem6.4}
Let  $D$ be a  bounded  domain that  satisfies an uniform exterior
cone condition, and let   $f$ be a  bounded  continuous function
on ${\bb R}^d$. If $(x_n,\alpha_n)$ converges  to $(x,\alpha)$ in
$D\times (0,2)$, then for each $T>0$,
\[
E^{\alpha_n}_{x_n}[\,f(X_{T\wedge\tau_D})\,] \to
E^{\alpha}_x[\,f(X_{T\wedge \tau_D})\,].
\]
\end{lemma}

\begin{proof}
Let
\[
{\cl C}_D(T) = {\cl C}_D \cap \{X(\tau^-_D) \in D\}\cap
\{\tau_D\ne T\}\cap \{\lim_{s\to T} X_s = X_T\}.
\]
Recall that the symmetric $\alpha$-stable process has no fixed
discontinuities. Then  by the eigenfunction expansion of
$P^{\alpha}_x(\tau_D>t)$,
\[
P^{\alpha}_x(\tau_D\ne T, \lim_{s\to T} X_s = X_T) = 1.
\]
Thanks to Lemma \ref{lem6.3},
\[
P^{\alpha}_x ({\cl C}_D(T)) = 1.
\]
If we can show
\[
\omega\in {\cl C}_D(T) \to \omega(T\wedge \tau_D(\omega))
\]
is continuous, then by an extension of the continuous mapping
theorem (Theorem 5.1 in \cite{Billingsley}),  the desired
conclusion will follow.

Let $\omega\in {\cl C}_D(T)$ and suppose $\omega_n$ converges to
$\omega$ in ${\bb D}[0,\infty)$. Define
\begin{align*}
t_n &= T\wedge \tau_D(\omega_n)\\ t &= T\wedge\tau_D(\omega).
\end{align*}
By Lemma \ref{lem6.2}, $\lim_{ n \to \infty}t_n =t$. Applying
Proposition 6.5 (a)  on page 125 in \cite{Ethier},
\[
\lim_{n \to \infty}|\omega_n(t_n)-\omega(t)| \wedge
|\omega_n(t_n)-\omega(t^-)|\,=\,  0.
\]

If $|\omega_n(t_n)-\omega(t)|$ converges to $0$, then clearly
\[
\lim_{n \to \infty}\omega_n(T\wedge \tau_D(\omega_n))
=\omega(T\wedge \tau_D(\omega)),
\]
as desired.

On the other hand, if
\begin{equation}\label{eq6.9}
\lim_{n \to \infty} |\omega_n(t_n)-\omega(t^-)|= 0,
\end{equation}
then we distinguish two cases:\ $T>\tau_D(\omega)$ and
$T<\tau_D(\omega)$ (recall $\omega\in {\cl C}_D(T)$ implies
$\tau_D(\omega)\ne T$).

Let us first assume $T>\tau_D(\omega)$. Since $\tau_D(\omega_n)$
converges to $\tau_D(\omega)$ by Lemma \ref{lem6.2},
$t_n=\tau_D(\omega_n)$ for large $n$. Hence by \eqref{eq6.9}
\[
\lim_{n\to\infty} \omega_n(\tau_D(\omega_n)) = \lim_{n\to\infty}
\omega_n(t_n) = \omega(t^-) = \omega(\tau_D(\omega)^-).
\]
Notice if $\lim\limits_{n\to\infty} \omega_n(\tau_D(\omega_n)) =
y$ exists, then $y\in D^c$. But then $\omega\in {\cl C}_D(T)$
implies  $y = \omega(\tau_D(\omega)^-)\in D$; contradiction. Thus
$T>\tau_D(\omega)$ is not possible.

Finally, if $T<\tau_D(\omega)$,  then by Lemma \ref{lem6.2}
\[
T<\tau_D(\omega_n),
\]
for $n$ large. Since $\omega\in {\cl C}_D(T)$, \eqref{eq6.9}
becomes
\[
\omega_n(T)\to \omega(T^-) = \omega(T).
\]
We  conclude
\begin{align*}
\lim_{n\to\infty} \omega_n(T\wedge \tau_D(\omega_n)) &=
\lim_{n\to\infty} \omega_n(T)\\ &= \omega(T^-)\\ &= \omega(T)\\ &=
\omega(T\wedge \tau_D(\omega)).
\end{align*}

In any event, we get the desired continuity.
\end{proof}

\section{Proof of Theorem \ref{thmb}}

\indent

We will need the following lemma; it is formula (2.7) in
\cite{Banuelos}. Though the authors  do not mention the statement
concerning continuity in $\alpha$, it is possible to trace back
through the literature they cite to see the statement holds.

\begin{lemma}\label{lem7.1}
If $D\subseteq {\bb R}^d$ is a bounded Lipschitz domain, then for
some positive continuous functions $C(\alpha)$ and
$\beta(\alpha)$,
\[
E^{\alpha}_x(\tau_D) \,\le\, C(\alpha)\,
\delta^{\beta(\alpha)}_D(x),\qquad \text{ for all } x\in D.\qquad
\]
\end{lemma}

The next result immediately follows.

\begin{corollary}\label{cor7.2}
Given a bounded Lipschitz domain $D$ and compact $K\times [a,b]
\subseteq \overline{D}\times(0,2)$,
\[
\sup\left\{\,E^{\alpha}_x[\tau_D]\colon \ (x,\alpha) \in K\times
[a,b]\;\right\} < \infty. \qquad\]
\end{corollary}

Corollary \ref{cor7.2} will allow us to get equicontinuity of the
eigenfunctions near $\partial D$. For the interior of $D$ we need
the following Krylov--Safanov type of theorem. Let
\[
G^{\alpha}_0g(x) = E^{\alpha}_x \left[\,\int^{\tau_D}_0
g(X_t)dt\,\right]
\]
be the 0-resolvent of the killed symmetric $\alpha$-stable process
in $D$.

\begin{lemma}\label{lem7.3}
Suppose $g$ is bounded with support in $\ovl D$. Then for each
$x\in D$ there exist positive continuous functions  $C(\alpha)$
and $\beta(\alpha)$, independent of $g$, such that  for all $
y\in$ D
\[
|\,G^{\alpha}_0g(x) - G^{\alpha}_0g(y)\,| \le
C(\alpha)\,\left[\,\sup|G^{\alpha}g| + \sup|g|\,\right]\,
|x-y|^{\beta(\alpha)}.
\]
\end{lemma}

\begin{proof}
This theorem is essentially due to Bass and Levin (see their
Proposition 4.2 on page 387). While they consider the 0-resolvent
\[
S_0 g(x) = E_x\left[\,\int^\infty_0 g(X_t)dt\,\right],
\]
their proof also works for the killed resolvent because their
crucial formula
\[
S_0 g(y) = E_y\left[\,\int_0^{\tau_{B(x,r)}} g(X_t)dt\,\right] +
E_y \left[\,S_0 g(X_{\tau_{B(x,r)}})\,\right]
\]
holds when $S_0$ is replaced by $G^\alpha_0$ and $E_y$ is replaced
by $E^{\alpha}_y$, where $r>0$ is such that $B(x,r) \subset D$.
Since we are restricted to $D$ instead of ${\bb R}^d$, the numbers
$C(\alpha)$ and $\beta(\alpha)$ depend on $x$, in contrast to the
case treated by Bass and Levin. Moreover, it is a simple matter to
go through their proof and see the numbers $C(\alpha)$ and
$\beta(\alpha)$ can be chosen to depend continuously on $\alpha$.
\end{proof}

\begin{corollary}\label{cor7.4}
Assume $D$ is bounded and Lipschitz. Then for each $x\in D$ and
$[a,b]\subseteq (0,2)$, there exist positive $C$ and $r$ such that
\[
|\,G^{\alpha}_0g(x) - G^{\alpha}_0g(y)\,| \le \,C
\,|x-y|^r\,\sup|g|
\]
for all $y\in D$, $\alpha \in [a,b]$ and bounded $g$ with support
in $\ovl D$.
\end{corollary}

\begin{proof}
By Corollary \ref{cor7.2}
\begin{align*}
\sup|G^{\alpha}g| &\le \sup|g|\cdot \sup_x E^{\alpha}_x(\tau_D)\\
&\le \sup|g|\cdot C
\end{align*}
where $C$ is independent of $\alpha\in [a,b]$ and $g$. The  result
follows from this and the continuity of $C(\alpha)$ and
$\beta(\alpha)$ from Theorem 7.3.
\end{proof}

At last we can prove Theorem \ref{thmb}.  It is well known that
\begin{equation}\label{eq7.1}
0\le p^{\alpha}_D(t,x,y) \le p^{\alpha}(t,x,y).
\end{equation}
Moreover,
\begin{equation}\label{eq7.2}
p^{\alpha}(t,x,y) \le C(\alpha) t^{-d/\alpha}
\end{equation}
where $C(\alpha)$ is continuous in $\alpha$ (see (2.1) in
\cite{Blumenthal}).

Let $\{\alpha_m\}_{m=1}^\infty$ be a  sequence converging to
$\alpha$. Recall that for all $ \beta \in (0,2)$,
$\{\varphi^{\beta}_m\}_{m=1}^\infty$ is an orthonormal set. Then
thanks to the symmetry of the heat kernel and (\ref{eq7.2})
\begin{align*}
\varphi^{\beta}_i(x) &= e^{\lambda^{\beta}_i t}\intl_D
p^{\beta}_D(t,x,y) \varphi^{\beta}_i(y)dy\\ &\le
e^{\lambda^{\beta}_it} \sqrt{\intl_D
\,\left[\,p^{\beta}_D(t,x,y)\,\right]^2dy}\\ &=
e^{\lambda^{\beta}_it} \sqrt{p^{\beta}_D(2t,x,x)}\\ &\le
e^{\lambda^{\beta}_it} \sqrt{\frac{C(\beta)}{(2t)^{d/\beta}}}.
\end{align*}
In particular, taking $t=1$ and using Theorem 1.2,
\begin{equation}\label{eq7.3}
\sup_{ \ x\in D, m\ge 1} \varphi^{\alpha_m}_i(x)\, \le\,
\sup_{m\ge 1}\,
e^{\lambda^{\alpha_m}_i}\sqrt{\frac{C(\alpha_m)}{2^{d/\alpha_m}}}
< \infty.
\end{equation}
Thus for  each $i\geq 1$, the sequence
$\left\{\varphi^{\alpha_m}_i\right\}_{ m= 1}^\infty$ is uniformly
bounded.  Next we show the sequence
$\left\{\varphi^{\alpha_m}_i\right\}_{ m= 1}^\infty$ is pointwise
equicontinuous on $\ovl D$. Indeed, since
\begin{equation}\label{eq7.4}
\varphi^{\beta}_i = \lambda^{\beta}_i G^{\beta}_0
\varphi^{\beta}_i,
\end{equation} Corollary \ref{cor7.4} implies that for each $x\in D$ there exist $C$ and $r$
such that
\begin{align*}
|\varphi^{\alpha_m}_i(x) - \varphi^{\alpha_m}_i(y)| &=
\,\lambda^{\alpha_m}_i\, |G^{\alpha_m}_0 \varphi^{\alpha_m}_i(x) -
G^{\alpha_m}_0 \varphi^{\alpha_m}_i(y)|\\ &\le \,C\,\left[\;
\sup_{m\ge 1} \,\left[\,\mu_i\,\right]^{\alpha_m/2}\, \right]\,
\,\left[\; \sup_{ u \in D,m\geq 1}
|\varphi^{\alpha_m}_i(u)|\,\right]\, |x-y|^r
\end{align*}
for all $m \ge 1$ and $y\in D$. Thanks to  \eqref{eq7.3} we get
the desired equicontinuity for $x\in D$.

As for $x\in\partial D$, first notice  \eqref{eq7.4} and Lemma
\ref{lem7.1} imply there are $r$ and $C$ independent of $m$ such
that for each $z\in D$
\begin{align*}
|\varphi^{\alpha_m}_i(z)| &\le \,\left[\sup_{m\ge 1}
\,\left[\,\mu_i\,\right]^{\alpha_m/2}\right]\, \left[\sup_{y \in
D,m\ge 1} |\varphi^{\alpha_m}_i(y)|\,\right]\,
E^{\alpha_m}_z(\,\tau_D\,)\\ &\le
C\,\left[\,\delta_D(z)\,\right]^r.
\end{align*}
Thus $\varphi^{\alpha_m}_i$ is continuous on $\ovl D$ with
boundary value 0. Hence if $x\in \partial D$ then
\begin{align*}
|\varphi^{\alpha_m}_i(x) - \varphi^{\alpha_m}_i(y)| &=
|\varphi^{\alpha_m}_i(y)|\\ &\le
C\,\left[\,\delta_D(y)\,\right]^r\\ &\le C|x-y|^r.
\end{align*}
By Ascoli's Theorem, the sequence
$\left\{\varphi^{\alpha_m}_i\right\}_{ m= 1}^\infty$  is
precompact in $C(\ovl D)$.

Next assume $\{\lambda^{\alpha_m}_i\}_{m=1}^{\infty}$ converges to
$\lambda$. We show any limit point $\varphi$ of the sequence
$\{\varphi^{\alpha_m}\}_{m=1}^{\infty}$ is an eigenfunction of
$H_\alpha$ and the corresponding eigenvalue is $\lambda$. Choose a
subsequence $\beta_r = \alpha_{m_r}$ such that, as $r\to\infty$,
$\varphi^{\beta_r}_i$ converges uniformly to $\varphi$ on $\ovl
D$. Since $\varphi^{\beta_r}_i$ and $\varphi$ are 0 on $\partial
D$, we can extend them to all of ${\bb R}^d$ by taking them to be
0 outside $D$. Then
\begin{align*}
E^{\beta_r}_x
\left[\,\varphi^{\beta_r}_i(X_{t\wedge\tau_D})\,\right] &=
E^{\beta_r}_x \left[\,\varphi^{\beta_r}_i(X_t)\;
I_{\tau_D>t}\,\right]\\
E^{\beta_r}_x\left[\,\varphi(X_{t\wedge\tau_D})\,\right] &=
E^{\beta_r}_x\left[\,\varphi(X_t)\; I_{\tau_D>t}\,\right],
\end{align*}
and $\varphi^{\beta_r}_i$ converges to $\varphi$ uniformly on
${\bb R}^d$. Thus we have
\begin{eqnarray}\nonumber
e^{-\lambda^{\beta_r}_it} \varphi^{\beta_r}_i(x) &= &\intl_D
p^{\beta_r}_D(t,x,y) \varphi^{\beta_r}_i(y)dy\\ \label{eq7.15} &=
&E^{\beta_r}_x \left[\,\varphi^{\beta_r}_i(X_t)\;
I_{\tau_D>t}\,\right]\\ \nonumber &=&
E^{\beta_r}_x\left[\,\varphi^{\beta_r}_i
(X_{t\wedge\tau_D})\,\right]\\ \nonumber &=&  E^{\beta_r}_x
\left[\,\varphi^{\beta_r}_i(X_{t\wedge\tau_D}) -
\varphi(X_{t\wedge\tau_D})\,\right] +
E^{\beta_r}_x\left[\,\varphi(X_{t\wedge\tau_D})\,\right]
\end{eqnarray}

Lemma \ref{lem6.4} and the uniform convergence of
$\varphi^{\alpha_m}_i$  to $\varphi$  imply

\[ \lim_{ r \to \infty} E^{\beta_r}_x \left[\,\varphi^{\beta_r}_i(X_{t\wedge\tau_D}) -
\varphi(X_{t\wedge\tau_D})\,\right] +
E^{\beta_r}_x\left[\,\varphi(X_{t\wedge\tau_D})\right]\,=
E^{\alpha}_x\left[\,\varphi(X_{t\wedge\tau_D})\,\right].\]

Since the left hand side of (\ref{eq7.15}) converges to
$e^{-\lambda t}\varphi(x)$, we conclude
\begin{align*}
e^{-\lambda  t} \varphi(x) &= E^{\alpha}_x
\left[\,\varphi(X_{t\wedge \tau_D})\,\right]\\ &= E^{\alpha}_x
\left[\,\varphi(X_t) \;I_{\tau_D>t}\,\right]\\ &= \intl_D
p^{\alpha}_D (t,x,y) \varphi(y)dy.
\end{align*}
Hence $\varphi$ is an eigenfunction of $H_\alpha$, and the
corresponding eigenvalue is $\lambda$. $\hfill\square$

\end{document}